%% file: article.tex
\documentclass{birkmult}

\usepackage{graphicx}

\newcommand{\eg}{{\it e.g.\/}}
\newcommand{\ie}{{\it i.e.\/}}
\newcommand{\R}{{\mathbb R}}
 
\newtheorem{thm}{Theorem}[section]

\theoremstyle{definition}
\newtheorem{defn}[thm]{Definition}
\theoremstyle{remark}

\newtheorem{ex}{Example}
\numberwithin{equation}{section}

\begin{document}

\title{Thinning Out Redundant Empirical Data}
\author{John Abbott}
\address{Dip. di Matematica, Universit\`a di Genova, via Dodecaneso 35, 16146 Genova, Italy}
\email{abbott@dima.unige.it}

\author{Claudia Fassino}
\address{Dip. di Matematica, Universit\`a di Genova, via Dodecaneso 35, 16146 Genova, Italy}
\email{fassino@dima.unige.it}

\author{Maria-Laura Torrente}
\address{Scuola Normale Superiore, piazza dei Cavalieri 7, 56126 Pisa, Italy}
\email{m.torrente@sns.it}
\date{}

\begin{abstract} Given a set $\mathbb X$ of ``empirical" points, whose coordinates are perturbed by errors,
we analyze whether it contains redundant information, that is whether some of its elements could be represented by a single equivalent point. If this is  the case, the empirical information associated to $\mathbb X$ could be described by fewer points, chosen in a suitable way.
We present two different methods to reduce the cardinality of $\mathbb X$ which compute a new  set of points equivalent to the original one, that is representing the same empirical information. Though our algorithms use some basic notions of Cluster Analysis they are specifically designed for ``thinning out" redundant data.
We include some experimental results which illustrate the practical effectiveness of our methods.
\end{abstract}

\maketitle

\input{introduction}

\input{equivalence}

\input{algorithms}
\input{clu_analysis_new}

\input{test}

\input{conclusions}

\bigskip
{\bf Acknowledgments}
Part of this work was conducted during the Special Semester on Gr\"obner Bases, February 1 to July 31, 2006, organized by the RICAM Institute (Austrian academy of Sciences) and the RISC Institute (Johannes Kepler University) in Linz, Austria.
We thank Lorenzo Robbiano for helpful discussions and suggestions on the subject of this paper.

\end{document}

%% file: introduction.tex
\section{Introduction}

Often numerical data in scientific computing arise from real-world
measurements, and so are perturbed by noise, uncertainty and approximation.
A common technique to counter this phenomenon is to make ``excessively
many'' measurements, and as a consequence the resulting body of empirical data
appears as a ``redundant'' set carrying relatively little information
compared to its cardinality.  Our aim is to reduce this redundancy by
replacing subsets of close values, which we regard as repeat measurements,
by a single representative value.

We view an empirical point $(p,\varepsilon)$ as a ``cloud'' of data which differ from $p$  by
less than the tolerance $\varepsilon$. 
If the intersection of different clouds is ``sufficiently" large, we can replace them by a single empirical point carrying essentially the same empirical information.
We illustrate this intuitive idea in the following example where an initial set of~$12$ points is ``thinned out" to an equivalent set of~$4$ points.
\begin{ex} \label{ex11}
Given the set $\mathbb X$ of $12$ points in $\R^2$
\begin{eqnarray*}
\mathbb X&=&\{(-1,-1),\;(0,-1),\;(1,-1),\;(-1,0),\;(0,0),\;(1,0),\\
& & \phantom{\{}(-1,1),\;(0,1),\;(1,1),\;(5,-2.9),\;(5,0),\;(5,2.9) \}
\end{eqnarray*}
we suppose that each coordinate is perturbed by an error less than $1.43$.
\begin{figure}
\begin{center}
\includegraphics[bb=14 14 432 349, width=4.5cm]{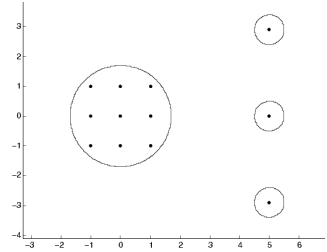}
\caption{Appropriate partition of $\mathbb X$}\label{fig1}
\end{center}
\end{figure}

In this situation, the first nine points most likely derive from measurements of the same
quantity; therefore it is quite reasonable (and appropriate) to 
collapse them onto a single candidate, for example the point $(0,\;0)$.
In contrast, since the last three points are well separated, they should not 
be collapsed.
This partition, shown in Figure~\ref{fig1},  is found by our algorithms, as reported in Examples \ref{ex31} and \ref{ex32}. 
\end{ex}

Based on the idea of clustering together empirical points which could derive from different
measurements of the same datum, we have designed  two algorithms which take a large set  of redundant data
and produce a smaller set of ``equivalent'' empirical points.  Typically the smaller set contains far fewer elements than the original one, with obvious consequent gains both in computational speed and in memory resources used in subsequent processing of the data. 

\smallskip
This paper is organized as follows. In Section $2$ we  introduce 
the concepts and tools useful to our work, focussing our attention on the idea of ``collapsable sets'' of empirical points.
Section $3$ describes the Agglomerative and the Divisive Algorithms to thin out sets of empirical points
while preserving the overall geometrical structure.
The relationship with the theory of Cluster Analysis
is discussed in Section~$4$.  In Section~$5$ we present some numerical examples to illustrate the behaviour of our algorithms on different geometrical configurations of points. The conclusions are summarized in Section 6.

%% file: equivalence.tex
\section{Basic Definitions and Notation}

This section recalls the definitions and tools used later in the paper.

We suppose that the points belong to the space  $\R^n$, $n \ge 1$, and we use
 the norm $\| \cdot \|_2$.  Further, given an $n\times n$ positive diagonal matrix
 $E$, we shall also use the weighted norm  $\| \cdot \|_{E,2}$ as defined in~\cite{DBA}.  For completeness, we recall here their definitions:
 \begin{eqnarray*}
\|v\|_2:= \sqrt{\sum_{j=1}^n v_j^2}\;\;\;\;{\rm and}\;\;\;\; \|v\|_{E,2}:=\|Ev\|_2
\end{eqnarray*}
Later on the index $2$ will be omitted for simplicity of notation.

\smallskip
Intuitively an empirical point, representing  real-world measurements, is a point $p$ of $\R^n$ whose coordinates are affected by noise, while  only an error estimation on them is known.
We suppose we know (for every $1 \leq i \leq n$) an estimate 
$\varepsilon_i \in \R^+$ of the error in the $i$-th component 
of  $p$, so that each point $r$ which differs from $p$ 
componentwise by less than $\varepsilon_i$ can be considered equivalent to 
$p$ from a numerical point of view. We can formalize this idea by means 
of the definition of empirical point, introduced by Stetter in~\cite{St}.

\smallskip
\begin{defn}\label{empirPts} 
Let $p \in \R^n$ be a point and let 
$\varepsilon=(\varepsilon_1,\dots, \varepsilon_n)$ with each $\varepsilon_i \in 
\R^+$, be the vector of the componentwise  estimated errors.
An {\bf empirical point} $p^{\varepsilon}$ is the pair $(p,\varepsilon)$, where
we call $p$ the {\bf specified value} and $\varepsilon$ the {\bf tolerance}.
\end{defn}

\smallskip
In this paper we shall consider sets of empirical points all having the same 
fixed tolerance $\varepsilon$.  This is a natural assumption if the 
points derive from real-world data measured with the same accuracy.
Additionally, this hypothesis simplifies the theoretical study.

From now on we denote by $\varepsilon=(\varepsilon_1,\dots,\varepsilon_n)$ 
with each $\varepsilon_i \in \R^+$, the fixed tolerance.  So given any $p \in \R^n$, we 
write $p^\varepsilon$ to mean the corresponding empirical point having $p$ as
specified value and $\varepsilon$ as tolerance.
We denote by $\mathbb{X}^{\varepsilon} = \{p_1^{\varepsilon},\ldots,
p_s^{\varepsilon}\}$ a set of empirical points each having the tolerance
$\varepsilon$ and by $\mathbb{X}=\{p_1,\ldots,p_s\}$ the set of the specified
values associated to $\mathbb{X}^{\varepsilon}$.
We define the diagonal matrix $E=diag(1/\varepsilon_1,\dots,1/\varepsilon_n)$ 
and shall use the $E$-weighted norm on $\R^n$ in order to ``normalize" 
the distance between points w.r.t.~the tolerance $\varepsilon$.

An empirical point $p^{\varepsilon}$ naturally defines the following set:
$$
N(p^\varepsilon) = \{r \in \R^n \;:\; \|p-r\|_{E} \le 1\} $$
Each element in $N(p^\varepsilon)$ can be
obtained by perturbing the coordinates of the specified value $p$ by amounts less than the tolerance; for this reason we can say that the points of $N(p^\varepsilon)$ represent the same
empirical information as $p$. Analogously, each element of $\cap_{p\in \mathbb X} N(p^\varepsilon)$, if this   intersection is not empty, represents the same empirical information as the whole set $\mathbb X$.
Although the choice of an element in this intersection is quite free, we decide to represent a set of ``close" points with their centroid. The following definitions are introduced in order to formalize this idea. 

\begin{defn}\label{def22}
The set of empirical points $\mathbb{X}^\varepsilon = \{p_1^\varepsilon,\ldots,p_s^\varepsilon\}$ is {\bf collapsable} if
\begin{eqnarray}\label{def_collaps}
\|p_i-q\|_E \le 1 \quad \forall i=1,\dots,s
\end{eqnarray}
where $q= \frac{1}{s}\sum_{i=1}^s p_i$ is the centroid of $\mathbb X$.
\end{defn}

If  $\mathbb X^\varepsilon$ is collapsable, the centroid $q$ of $\mathbb X$ 
belongs to each of the sets $N(p_i^\varepsilon)$; so the empirical point 
$q^\varepsilon$ is numerically equivalent to every point in 
$\mathbb X^\varepsilon$. 
We formalize this idea as follows.
\begin{defn}\label{def23}
The {\bf empirical centroid} of a set $\mathbb X^\varepsilon$ is the 
empirical point $q^\varepsilon$ where~$q$ is the centroid of the
set $\mathbb X$.  If $\mathbb X^{\varepsilon}$ is a collapsable set, its 
empirical centroid is called its {\bf valid representative}.
\end{defn}

If a set of empirical points contains a collapsable subset, it contains
some redundancy, {\it i.e.}~it carries relatively little empirical
information compared to number of points in it.  The methods presented in
this paper are designed to ``thin out" such sets by finding a smaller set
of empirical points with much lower redundancy which still contain
essentially the same empirical information.

%% file: algorithms.tex
\section{Algorithms }

In this section we describe two algorithms that, given a set 
$\mathbb X^\varepsilon$ of empirical points, compute a partition 
$\mathcal L^\varepsilon =\{L_1^\varepsilon,\dots,L_k^\varepsilon\}$ of it, consisting of 
non-empty collapsable sets, and a set $\mathbb Y^\varepsilon=
\{q_1^\varepsilon,\dots,q_k^\varepsilon\}$ where each $q_i^\varepsilon$ is the 
valid representative of $L_i^\varepsilon$.
Our algorithms differ in the strategies for building the partitions:
\begin{enumerate}
\item the {\bf Agglomerative Algorithm} initially puts each point of 
$\mathbb X^\varepsilon$ into a different subset and then iteratively unifies 
pairs of subsets into a larger collapsable set;
\item the {\bf Divisive Algorithm} initially puts all the points of 
$\mathbb X^\varepsilon$ into a single subset and then iteratively splits off 
the remotest outlier and ``evens up'' the new partition.
\end{enumerate}

\subsection{The Agglomerative Algorithm}
The Agglomerative Algorithm (AA) implements a unifying method. The sets in the partition are determined by an iterative process.  Initially each set contains a single original empirical point, then iteratively the two 
closest sets are unified provided their union is collapsable. This method is quite fast when the input points are well separated w.r.t.~the tolerance, since a small number of set unifications is required.  

\begin{thm}(The Agglomerative Algorithm)\\
Let $\mathbb X^\varepsilon=\{p_1^\varepsilon,\dots,p_s^\varepsilon\} $ be a set of
empirical points, with each $p_i \in \R^n$ and a common tolerance 
$\varepsilon=(\varepsilon_1,\dots,\varepsilon_n)$. 
 Let $\|\cdot\|_{E}$ be the weighted norm on $\R^n$ 
w.r.t. $E~=~diag(1/\varepsilon_1,\dots,1/\varepsilon_n)$.  Consider the following sequence of instructions.  
\begin{description}
\item[AA1] Start with the subset list $\mathcal L=[L_1,\dots,L_s]$ where each
$L_i=\{p_i\}$, and the list ${\mathbb Y}=[q_1,\dots,q_s]$ of the centroids
of the $L_i$. 
\item[AA2]  Compute the symmetric matrix $M=(m_{ij})$ such that 
$m_{ij}=\|q_i-q_j\|_E$ for each $q_i,q_j \in \mathbb Y$.
\item[AA3] If $|\mathbb Y| = 1$ or $ \min \{m_{ij} : i<j\}>2$ 
then return the lists $\mathcal L$ and $\mathbb Y$ and stop.
\item [AA4] Choose $\hat {\i}$, $\hat {\j}$ s.t. $m_{\hat{\i} \hat{\j}}=
\min \{m_{ij} : i<j\}$ and compute the centroid $q$ of $L_{\hat{\i}}\cup L_{\hat{\j}}$
$$ q =\frac{ |L_{\hat {\i}}| q_{\hat {\i}} + |L_{\hat {\j}}| q_{\hat {\j}}}{ |L_{\hat {\i}}| +  |L_{\hat {\j}}| }$$
\item [AA5] If  $\|p-q\|_E \le  1$ for every $p \in L_{\hat {\i}}\cup L_{\hat {\j}}$ then in $\mathcal L$  replace $L_{\hat {\i}}$ by  
$L_{\hat {\i}} \cup L_{\hat {\j}}$ and remove $L_{\hat {\j}}$. Similarly, 
in $\mathbb Y$ replace $q_{\hat {\i}}$ by $q$ and remove $q_{\hat {\j}}$ 
and then go to step~AA2.
Otherwise put $m_{\hat{\i} \hat{\j}}=\infty$ (any value greater than $2$
will do) and go to step~AA3.
\end{description}
This algorithm computes a pair $(\mathcal L,\;\mathbb Y)$ such that:
\begin{itemize} 
\item $\{ L_i^\varepsilon : L_i \in \mathcal L\}$ is a partition of $\mathbb X^\varepsilon$
into collapsable sets such that no pair can be unified into a collapsable set;
\item for each $q_i \in \mathbb Y$ the empirical point $q_i^\varepsilon$ is
the valid representative of $L_i^\varepsilon$. 
\end{itemize}
\end{thm}
\begin{proof} First we prove finiteness.  Step~AA2 is performed only
  finitely many times and so a finite number of matrices $M$ is computed.
  In fact, after the first computation of $M$, this step is performed only
  when the algorithm removes an element from $\mathbb Y$, \ie~at most $s-1$
  times.  Now, also step~AA4 is performed only finitely many times on the
  same matrix $M$, because it is performed only when the minimal element
  $m_{\hat{\i} \hat {\j}}$ of the matrix $M$ is less than or equal to $2$
  and then either two subsets are unified or $m_{\hat{\i} \hat{\j}}$ is
  replaced by $\infty$, but this can happen at most $s^2/2$ times.

Next we show correctness.  First, note that the elements of $\mathcal L$ 
define a partition of $\mathbb X$.  In fact, in step~AA1 we set $\mathcal L=
[\{p_1\},\dots,\{p_s\}]$; the only place where $\mathcal L$ changes is in Step~AA5 when
we unite two of its elements, and so a new partition of 
$\mathbb X$ is obtained.  Obviously $\mathcal L^\varepsilon$ is also a partition 
of $\mathbb X^\varepsilon$.

For each $L_i \in \mathcal L$, the corresponding 
empirical set $L_i^\varepsilon$ is collapsable.  This is clearly true in step~AA1.
Step~AA5 unites two elements of $\mathcal L$ only if their union is collapsable:
step~AA4 computes the centroid 
$q$ of $L_i \cup L_j$ and step~AA5 tests condition~(\ref{def_collaps}) for 
each point in $L_i\cup L_j$.

Now we prove that upon termination the union of any pair of elements of
$\mathcal L$ is not collapsable.  If the algorithm stops because $\mathbb
Y$ (and $\mathcal L$ too) contains a single element, the conclusion is
trivial.  Otherwise, the algorithm ends because $m_{ij} > 2$ for all $i<j$.
We observe that the elements $m_{ij}$ of the final matrix $M$ are such that
either $m_{ij}=\|q_i-q_j\|_E$ or $m_{ij}=\infty$ but $\|q_i-q_j\|_E \le 2$.
The case where $m_{ij}=\infty$ is trivial: an entry in $M$ can become $\infty$
only in step~AA5 after having verified that $L^\varepsilon_i \cup
L^\varepsilon_j$ is not collapsable.  In the case where $m_{ij}$ is finite we
show that the union of $L^\varepsilon_i$, $L^\varepsilon_j$ is a not collapsable
set by contradiction.  We suppose that $\|p-q\|_E \le 1$ for each $p \in
L_i\cup L_j$, where $q$ is the centroid of $L_i \cup L_j$.  If $m=|L_i|$
and $n=|L_j|$, we have
\begin{eqnarray*}
& &\|q_i-q_j\|_E= \left \| \frac{1}{m} \left( \sum_{p \in L_i} p - m q \right) + \frac{1}{n}\left( nq -\sum_{p \in L_j} p\right ) \right \|_E\\
&=& \left \| \frac{1}{m} \sum_{p \in L_i} (p - q ) + \frac{1}{n} \sum_{p \in L_j}( q-p ) \right \|_E
 \le  \frac{1}{m} \sum_{p \in L_i} \|p - q\|_E+ \frac{1}{n} \sum_{p \in L_j}\| q-p \|_E\end{eqnarray*}
From the hypothesis, we deduce that $\|q_i-q_j\|_E \le 2$, a contradiction.

Finally, we can conclude the proof since,  by construction, each element 
$q_i \in \mathbb Y$ is the centroid of $L_i$ and $L_i^\varepsilon$ is 
collapsable, so the empirical centroid $q_i^\varepsilon$ is indeed the valid 
representative of $L_i^\varepsilon$. 
\end{proof}

\smallskip
Note that, in step~AA5, we must check the condition that $\|p-q\|_E \le 1$
for each $p \in L_{\hat {\i}}\cup L_{\hat {\j}}$.  In fact, if we check only the
condition $\|q_{\hat{\i}}- q_{\hat{\j}}\|_E\le 1$, there are pathological
examples where not collapsable sets are built in the final partition
(see Example~\ref{ex_zip}).

The algorithm as presented here can easily be improved from the
computational point of view: in step~AA2 it is not necessary to compute a
new matrix $M$ after uniting $L_{\hat {\i}}$ and $L_{\hat {\j}}$, but
suffices to remove the $\hat {\j}$-th column and to update the $\hat
{\i}$-th row.

\smallskip
In the following example we apply the Agglomerative Algorithm on the 
points of Example \ref{ex11} to show that the desired partition is 
obtained (see Figure~\ref{fig1}).
\begin{ex}\label{ex31}
Let $\mathbb X^\varepsilon=\{p_1^\varepsilon,\dots,p_{12}^\varepsilon\}$ be a set 
of empirical points with tolerance $\varepsilon=(1.43,\; 1.43)$, whose 
specified values coincide with the set $\mathbb X$ of Example \ref{ex11}:
\begin{eqnarray*}
\mathbb X&=&\{(-1,-1), \;(0,-1),\;(1,-1),\;(-1,0),\;(0,0),\;(1,0),\\
& & \phantom{\{}(-1,1),\;(0,1),\;(1,1),\;(5,-2.9),\;(5,0),\;(5,2.9) \}
\end{eqnarray*}
The AA computes, at each step, the following partitions, only clustering 
together the first nine points.
\begin{enumerate}
\item  $\mathcal L=\bigl\{ \{p_1\},  \{p_2\},  \{p_3\}, \{p_4\}, \{p_5\}, \{p_6\},  \{p_7\}, \{p_8\}, \{p_9\},  \{p_{10}\}, \{p_{11}\}, \{p_{12}\} \bigr\}$
\item  $\mathcal L=\bigl\{{\bf  \{p_1, p_2\}},  \{p_3\}, \{p_4\}, \{p_5\}, \{p_6\}, \{p_7\}, \{p_8\},  \{p_9\},  \{p_{10}\}, \{p_{11}\}, \{p_{12}\}\bigr\}$
\item  $\mathcal L=\bigl\{ {\bf \{p_1, p_2, p_4\}}, \{p_3\}, \{p_5\}, \{p_6\}, \{p_7\}, \{p_8\}, \{p_9\},  \{p_{10}\}, \{p_{11}\}, \{p_{12}\}
\bigr\}$
\item  $\mathcal L=\bigl\{ \{p_1, p_2, p_4\}, {\bf \{p_3, p_6\}}, \{p_5\}, \{p_7\}, \{p_8\}, \{p_9\},  \{p_{10}\}, \{p_{11}\}, \{p_{12}\}
\bigr\}$
\item  $\mathcal L=\bigl\{ \{p_1, p_2, p_4\}, \{p_3, p_6\}, {\bf\{p_5, p_8\}}, \{p_7\}, \{p_9\},  \{p_{10}\}, \{p_{11}\}, \{p_{12}\}
\bigr\}$
\item  $\mathcal L=\bigl\{{\bf \{p_1, p_2, p_4,p_5,p_8\}}, \{p_3, p_6\}, \{p_7\}, \{p_9\},  \{p_{10}\}, \{p_{11}\}, \{p_{12}\}\bigr\}$.
\item  $\mathcal L=\bigl\{ {\bf \{p_1, p_2, p_3, p_4,p_5,p_6, p_8\}}, \{p_7\}, \{p_9\},  \{p_{10}\}, \{p_{11}\}, \{p_{12}\}
\bigr\}$
\item  $\mathcal L=\bigl\{{\bf \{p_1, p_2, p_3, p_4, p_5, p_6, p_7, p_8\}}, \{p_9\},  \{p_{10}\}, \{p_{11}\}, \{p_{12}\}
\bigr\}$
\item  $\mathcal L=\bigl\{{\bf \{p_1, p_2, p_3, p_4, p_5, p_6, p_7, p_8, p_9\}},  \{p_{10}\}, \{p_{11}\}, \{p_{12}\}
\bigr\}$
\end{enumerate}
\end{ex}

\subsection{The Divisive Algorithm}

The Divisive Algorithm (DA) implements a ``subdivision" method.  The
sets in the partition are determined by an iterative process.  Initially
the partition consists of a single set containing all the points.  Then
iteratively DA seeks the original point farthest from the centroid of its
set.  If the distance between them is below the tolerance threshold then
the algorithm stops, because all original points are sufficiently well
represented by the centroids of their sets.  Otherwise it splits off the
worst represented original point into a new set initially containing just
itself.  Then DA proceeds with a redistribuition phase with the aim of
associating each original point to the current best representative subset
(locally) minimizing the total central sum of squares, defined as follows~\cite{RR}.

\begin{defn}
  Let $\mathbb{X}$ be a subset of $\R^n$ and let $q$ be its centroid. 
The {\bf central sum of squares} of $\mathbb X$ is defined to be:
$$
\sum_{p \in \mathbb{X}} \|p - q\|^2
$$
\end{defn}

\begin{defn}
  Let  $\mathcal L=\{L_1,\dots,L_k\}$ be a partition of the set $\mathbb{X}$.  The {\bf total
 central sum of squares} of the partition $\mathcal L$ is defined to be:
$$
I({\mathcal L}) = \sum_{j=1}^k I_j
$$
where $I_j$ is the central sum of squares of $L_j$.
\end{defn}

\smallskip
If $\mathbb X^\varepsilon$ contains large subsets of close empirical points,
DA turns out to be more efficient than AA, since a smaller number of
subdivisions is required.

\begin{thm}\label{pa}(The Divisive Algorithm)\\
Let $\mathbb X^\varepsilon=\{p_1^\varepsilon,\dots,p_s^\varepsilon\} $ be a set of empirical points, with each $p_i \in \R^n$ and a common tolerance 
$\varepsilon=(\varepsilon_1,\dots,\varepsilon_n)$.  Let $\|\cdot\|_{E}$ be the weighted norm on $\R^n$  w.r.t. $E=diag(1/\varepsilon_1,\dots,1/\varepsilon_n)$.  Consider the following sequence of instructions.  
\begin{description}
\item[DA1] Start with the list $\mathcal L = [L_1]$ where $L_1=\mathbb X$, and the centroid list~$\mathbb Y=[q_1]$ where~$q_1$ is the
 centroid of $L_1$.
\item[DA2] Let $\mathcal L=[L_1,\dots,L_r]$ and 
$\mathbb Y= [q_1,\dots,q_r]$,  the centroid list of the elements of 
$\mathcal L$. 
For each $p_i\in \mathbb X$ set $d_i=\| p_i - q_j\|_E$ where $L_j$ is the subset (of $\mathbb X$) to which 
$p_i$ belongs. Build the list $D = [d_1,\dots,d_s]$.
\item[DA3] If $\max(D) \le 1$ then return the lists 
$\mathcal L$ and $\mathbb Y$, and stop.
\item[DA4] Choose an index $\hat{\i}$ such that $d_{\hat{\i}}=\max(D)$, and compute  the index $\hat {\j}$ of the subset~$L_{\hat {\j}}$ to which $p_{\hat {\i}}$ belongs.
Remove $p_{\hat {\i}}$ from $L_{\hat {\j}}$ and compute the new centroid~$q_{\hat {\j}}$ of~$L_{\hat {\j}}$; 
append $L_{r+1}=\{p_{\hat {\i}}\}$ to 
$\mathcal L$ and $q_{r+1} =p_{\hat {\i}}$ to $\mathbb Y$.
\item[DA5] Compute the total central sum of squares $I(\mathcal L)$ of
the new partition $\mathcal L$.
\item[DA6]  For each $p \in \mathbb X$ and for each $L_k \in \mathcal L$, denote by
${\mathcal L}_{p,k}$ the partition  $\mathcal L$ but with $p$ moved into $L_k$.
Compute the total central sum of squares $I({\mathcal L}_{p,k})$.
\item[DA7] Choose a point $\hat{p} \in \mathbb X$ and an index 
$\hat{k}$ s.t.
$$
I({\mathcal L}_{\hat p,\hat k}) = \min \{I({\mathcal L}_{p,k}): p \in {\mathbb X},\; {L_k}\in\mathcal L \}
$$
\item[DA8] If $I({\mathcal L}_{\hat p,\hat k})  \ge I(\mathcal L)$ then go to~DA2.  Otherwise set
  $\mathcal L= {\mathcal L}_{\hat p,\hat k} $.  Compute the  centroids of the new partition $\mathcal L$. Go to~DA5.
\end{description}
This algorithm computes a pair $(\mathcal L,\;\mathbb Y)$ such that:
\begin{itemize} 
\item $\{ L_i^\varepsilon : L_i \in \mathcal L\}$ is a partition of $\mathbb X^\varepsilon$
into collapsable sets;
\item for each $q_i \in \mathbb Y$, the empirical point $q_i^\varepsilon$ is the valid representative of $L_i^\varepsilon$. 
\end{itemize}
\end{thm}
\begin{proof} 
  Later on we shall refer to the loop DA5--DA8 as ``the redistribution
  phase": points are moved from one subset to another in order to strictly decrease
  the total central sum of squares.  Note that in the redistribution phase the
  cardinality of~$\mathcal L$ does not change as the algorithm never
  eliminates any set in $\mathcal L$.  Indeed, if the singleton set $L_j=
  \{p\}$ belongs to $\mathcal L$, the point $p$ will not be moved to
  another set~$L_k \in \mathcal L$ leaving $L_j$ empty, since this new
  configuration cannot have smaller total central sum of squares: the combined
  central sum of squares of the sets~$L_j=\{p\}$ and $L_k$ is
$$I_j+I_k= 0+ \sum_{r\in L_k}\|r-q_k\|^2$$ 
where $q_k$ is the centroid of $L_k$, whereas the combined central sum of squares
of the new sets $L'_j=\emptyset$ and $L'_k=L_k \cup \{p\}$ is 
$$I'_j + I'_k = 0+ \left( \sum_{r\in L_k}\|r-q'_k\|^2 + \|p-q'_k\|^2 \right)$$ 
where $q'_k$ is the centroid of $L'_k=L_k \cup \{ p\}$.  And since $q_k$ 
is the centroid of $L_k$, we have $\sum_{r\in L_k}\|r-q'_k\|^2 \ge 
\sum_{r\in L_k}\|r- q_k\|^2$.  Consequently the new total central sum of squares cannot be smaller.

\smallskip
Now we prove finiteness.  The algorithm comprises two nested loops: the outer loop spanning
steps~DA2--DA8, and the redistribution phase (steps DA5--DA8).  The outer loop
cannot perform more than $s$ iterations because step~DA4 can be performed at
most $s$ times; anyway, after $s$ iterations the termination criterion in
step~DA3 will surely be satisfied as all the $d_i$ would be zero.

The redistribution loop will perform only finitely many iterations.  Each
iteration strictly reduces the total central sum of squares, and since
$\mathbb X$ is finite it has only finitely many partitions.
Consequently there are only finitely many possible values for the total
central sum of squares.

\smallskip
Next we show correctness.  The elements of $\mathcal L$ define a partition 
of $\mathbb X$.   This is trivially true in step~DA1.  The creation of a new
subset in step~DA4 clearly maintains the property.  The redistribution 
phase merely moves points between subsets (in step~DA8), so also preserves
the property.

The test in step~DA3 guarantees that upon completion of the algorithm each~$L_i \in \mathcal
L$ corresponds to a collapsable $L_i^\varepsilon$.  By construction, each
element $q_i \in \mathbb Y$ is the centroid of $L_i$.  Thus $q_i^\varepsilon$
is the valid representative of $L_i^\varepsilon$.
\end{proof}

\smallskip
In the following example we apply the Divisive Algorithm to the points of Example~\ref{ex11} to show that the desired partition is obtained (see Figure~\ref{fig1}). 

\begin{ex}\label{ex32}
Let $\mathbb X^\varepsilon=\{p_1^\varepsilon,\dots,p_{12}^\varepsilon\}$ be a set of empirical points with tolerance $\varepsilon=(1.43,\; 1.43)$, whose specified values coincide with the set $\mathbb X$ of Example \ref{ex11}:
\begin{eqnarray*}
\mathbb X&=&\{(-1,-1), \;(0,-1),\;(1,-1),\;(-1,0),\;(0,0),\;(1,0),\\
& &\phantom{\{}(-1,1),\;(0,1),\;(1,1),\;(5,-2.9),\;(5,0),\;(5,2.9) \}
\end{eqnarray*}
\goodbreak
The DA computes, at each step, after the redistribution phase, the 
following partitions.
\begin{enumerate}
\item  $\mathcal L=\bigl\{ \{p_1, p_2, p_3, p_4, p_5 , p_6, p_7, p_8, p_9, p_{10}, p_{11}, p_{12}\} \bigr\}$
\item  $\mathcal L=\bigl\{ \{{\bf p_1, p_2, p_3, p_4, p_5 , p_6, p_7, p_8, p_9}\},  \{{\bf p_{10}, p_{11}, p_{12}}\} \bigr\}$
\item  $\mathcal L=\bigl\{ \{p_1, p_2, p_3, p_4, p_5 , p_6, p_7, p_8, p_9\}, \{{\bf p_{10}}\},\{{\bf p_{11}, p_{12}}\} \bigr\}$
\item  $\mathcal L=\bigl\{ \{p_1, p_2, p_4, p_5,p_8,p_3, p_6,p_7,p_9\},  \{p_{10}\}, \{{\bf p_{11}}\}, \{{\bf p_{12}}\}\bigr\}$
\end{enumerate}
As mentioned before, DA performs fewer iterations than AA since several
input points are close together w.r.t. the tolerance. 
\end{ex}
\goodbreak

\subsection{A particularly quick method: the grid algorithm}\label{grid_sect}
We recall the $\infty$-norm and its corresponding $E$-weighted norm on $\R^n$, see~\cite{DBA}:
$$\|v\|_\infty=\max_{i=1\dots n} |v_i|\;\;\;\; {\rm and} \;\;\;\;
\|v\|_{E,\infty}=\|Ev\|_\infty
$$ 
where $E=diag(1/\varepsilon_1, \dots,1/\varepsilon_n)$, as before.

A particularly quick method for decreasing the cardinality of the set $\mathbb
X^\varepsilon $ can be designed using a regular grid, consisting of half-open balls of
radius~$1/2$ w.r.t.~the $E$-weighted norm $ \|\cdot\|_{E,\infty}$.  We
arbitrarily choose one ball to have the origin as its centre then
tessellate to cover the whole space.

This algorithm computes a partition of $\mathbb X^\varepsilon$ by gathering all
the empirical points whose specified values lye in the same ball into the same subset.
Suppose that one of these subsets comprises the empirical points $p_1^\varepsilon,\dots,p_m^\varepsilon$,
and let $q^\varepsilon$ be their empirical centroid, then $q^\varepsilon$ is a ``good" representative of each $p_i^\varepsilon$ because
$$\|p_i-q\|_{E,\infty}= \left \|p_i-\frac{1}{m}\sum_{j=1}^m p_j\right \|_{E,\infty} \le \frac{1}{m}\sum_{j=1}^m \|p_i-p_j\|_{E,\infty} < 1$$ 
However, in general such a subset is not collapsable, a notion defined in terms of the $2$-norm.

Note that, since the separations of the empirical points are ignored by
this method, unsatisfactory partitions can be obtained, \eg~close points
may happen to belong to different balls and so be assigned to different
subsets.  Nevertheless, this drawback is compensated by the speed and
simplicity of the method.  In particular, this grid method (with a smaller
radius) can be used to reduce the bulk of a very large body of data before
applying one of the more sophisticated but slower algorithms, AA or DA.
Another application of the grid method is to help choose the more suitable
algorithm between AA and DA by estimating the numbers of sets in the
partitions which would be produced.

%% file: clu_analysis_new.tex
\section{Relationship with Cluster Analysis}

The idea of analyzing a large body of empirical data and of partitioning
it into sets of ``similar values'' has been well studied in the theory of
Cluster Analysis (\eg~see~\cite{KM}).  The overall aim of Cluster Analysis
is to separate the original data into clusters where the members of each
cluster are much more similar to each other than to members of other clusters.
In contrast, our methods are more concerned with thinning out groups of
very close values while ignoring more distant points.
Below we show how Ward's ``classical" algorithm~\cite{RR}, an agglomerative
hierarchical method, and Li's more recent algorithm~\cite{Li}, a divisive hierarchical method, partition the empirical points of Example~\ref{ex11}.

\begin{ex}
  Let $\mathbb X^\varepsilon$ be the set of empirical points whose set
  of specified values is given in Example~\ref{ex11};
  similarly, let $\varepsilon=(1.43,1.43)$ as given there.  We recall
  that in Examples~\ref{ex31} and~\ref{ex32} both our algorithms AA
  and DA obtained the minimal partition into collapsable sets, as
  illustrated in Figure~\ref{fig1}.

Ward's and Li's algorithms do not obtain this minimal partition.
In fact, after~$8$ steps, Ward's algorithm puts the points $(5,-2.9)$ and
$(5,0)$ into the same cluster, while the first nine points of $\mathbb X$
still belong to different clusters.  Since this is an agglomerative method
no set of points is split during the computation, so Ward's
algorithm fails to recognise the collapsable set of nine points.  In a similar
vein, Li's algorithm goes astray at the third step: it divides the first
nine points of $\mathbb X$ into two subsets while the points $(5,-2.9)$ and
$(5,0)$ still belong to the same cluster.  Since this is a hierarchical divisive
method, once a set is split it can never be joined together again, so Li's algorithm
needlessly splits the collapsable set of nine points.
\end{ex}

\smallskip
Now we consider another method of Cluster Analysis, QT~Clustering~\cite{HKY},
because it has a number of similarities to our methods, especially AA.  QT~Clustering
computes a partition of the input data using a given limit on the diameter of the
clusters.  It works by building clusters according to their cardinality, while we are
primarily interested in the local geometrical separations of the input data.

\begin{ex}
  Let $\mathbb X^\varepsilon$ be a set of empirical points with tolerance
  $\varepsilon=(0.5)$ and with specified values
  $\mathbb X =\{ 0, \;0.05,\; 0.9,\;1,\;1.2 \} \subseteq \R$.
  Applying the QT Clustering algorithm with maximum cluster
  diameter equal to~$2\varepsilon$, we obtain the partition~$\bigl\{\{ 0,
    0.05,0.9,1\}, \; \{1.2\}\bigr\} $ where $\{ 0,
    0.05,0.9,1\}^\varepsilon$ is a not collapsable set.
  In contrast, if we apply AA or DA to $\mathbb X^\varepsilon$, we obtain the more
  balanced partition $\bigl\{ \{0, 0.05\},\; \{0.9,1,1.2\}\bigr\}$
  whose elements consist of specified values of collapsable sets. 
  We maintain that our partition is more
  plausible as a grouping of noisy data.
\end{ex}

%% file: test.tex
\section{Numerical Tests and Illustrative Examples}

In this section we present some numerical examples to show the
effectiveness and the potential of our techniques.  Both AA and DA
have been implemented using the C++ language, and are included in
CoCoALib~\cite{Co}.  All computations in the following examples have
been performed on an Intel Pentium~M735 processor (at~1.7~GHz) running
GNU/Linux and using the implementation in CoCoALib.

\bigskip
\begin{ex}{\bf Clouds of empirical points.}\\
  In this example we consider an empirical set $\mathbb X^\varepsilon$
  containing two well separated empirical points and three clusters,
  two big and one small.  Both AA and DA compute five valid
  representatives for $\mathbb X^\varepsilon$, but because the result
  comprises very few points DA is faster than AA.\\
Let $\mathbb X^\varepsilon$ be a set of empirical points, with tolerance  $\varepsilon=(20,20)$ and specified values
$\mathbb X=\cup_{i=1}^5 \mathbb X_i \subseteq \R^2$, where
\begin{itemize} 
\item[] $\mathbb X_1$ consists of $82$ points lying inside the disk of radius $10$ centered on $(0,0)$, 
\item[] $\mathbb X_2$ consists of  $64$  points lying inside the disk of radius $10$ centered on $(40,50)$, 
\item[] $\mathbb X_3=\{(49,0),(50,0),(50,1)\}$, $\mathbb X_4=\{ (9,41)\}$ and  $\mathbb X_5=\{ (-10,80)\}$.
\end{itemize}
Both AA and DA compute the ``intuitive" partition consisting of $5$ subsets $L_i=\mathbb X_i$ for $i=1,\ldots,5$, as shown in Figure \ref{clouds}.

\begin{figure}[h]
\begin{center}
\includegraphics[bb=14 14 428 320, width=6cm]{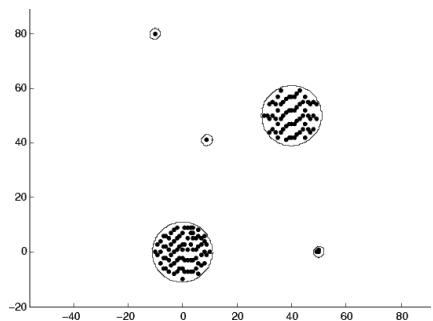}
\vspace{-.5cm}
\caption{Appropriate partition of $\mathbb X$}\label{clouds}
\end{center}
\end{figure}
\end{ex}

\begin{ex}{\bf Empirical points close to a circle.}\\
  In this example we compare the behaviour of AA and DA on a family of
  test cases, which comprises sets of empirical points with
  similar geometrical configurations but with differing ``densities".
  Let $\mathbb{X}_1 $, $\mathbb{X}_2 \subset \R^2$ be two sets of
  points lying close to the circle of radius $200$ and centered at the
  origin.  They contain $2504$ and $5032$ points, respectively.  For
  simplicity we choose a tolerance $\varepsilon=(\varepsilon_1,\varepsilon_2)$
  with $\varepsilon_1=\varepsilon_2$.  The numerical tests are performed by
  applying both AA and DA to the empirical sets $\mathbb{X}^\varepsilon_1
  $ and $\mathbb{X}_2^\varepsilon $ for various values of $\varepsilon$,
  viz.~$\varepsilon_1=2^k$ for $k=0,\dots,6$, since, for a fixed set of points simply increasing $\varepsilon$
  effectively increases the density of the points.
 
  In Table~\ref{tab1} we  present the results obtained processing $\mathbb X_1$ and $\mathbb
  X_2$ respectively.  The first column contains the value of the
  tolerance, the columns labeled with ``$\#$VR'' contain the number of
  the valid representatives computed by AA and DA respectively, while
  those labeled with ``Time" show the timings (in seconds) of each
  algorithm. The results show that DA runs quickly if $\varepsilon$ is large, that is
when the set of empirical points is dense enough, since only a few
splittings of the original set are needed.  On the other hand, when
the points are well separated, AA is preferable since the final
partition consists of a large number of sets.

\begin{table}[h]
\begin{tabular}{|r|r|r||r|r||r|r||r|r|}\cline{2-9} 
 \multicolumn{1}{c|}{}& \multicolumn{4}{c||}{\bf 2504 empirical points}& \multicolumn{4}{c|}{\bf 5032 empirical points}\\  
\cline{2-9}
 \multicolumn{1}{c|}{}& \multicolumn{2}{c||}{\bf AA} &  \multicolumn{2}{c||}{\bf DA}& \multicolumn{2}{c||}{\bf AA} &  \multicolumn{2}{c|}{\bf DA}\\  
\hline
$\varepsilon$ & $\#$VR & Time & $\#$VR & Time & $\#$VR & Time & $\#$VR & Time\\ \hline 
1 & 911 &  1 s & 727 & 293 s & 2096 & 6 s & 1460 & 2306 s\\
2 & 462 & 3 s & 347 & 184 s & 734 & 31 s & 587 & 1250 s\\
4 & 224 & 8 s & 173 & 114 s & 263 & 118 s & 185 & 577 s\\
8 & 108 & 18 s & 87 & 66 s & 121 & 317 s & 86 & 314 s\\  
16 & 56  & 50 s & 41 & 33 s &61  & 733 s & 41 & 166 s\\
32 & 29  & 117 s & 20 & 15 s & 28 & 1680 s & 21 & 79 s\\ 
64 & 13 & 2633 s & 10 & 6 s & 14 & 3695 s & 10 & 25 s \\
\hline
\end{tabular}
\vspace{0.1cm}
\caption{Points close to a circle }\label{tab1}
\vspace{-0.5cm}
\end{table}

  Figure~\ref{ring} shows a subset of $\mathbb X_1$ (the crosses) and its valid representatives (the~dots) w.r.t. the tolerance $\varepsilon=(16, 16)$.     

\begin{figure}[h]
\begin{center}
\includegraphics[bb=14 14 582 468, width=6.5cm]{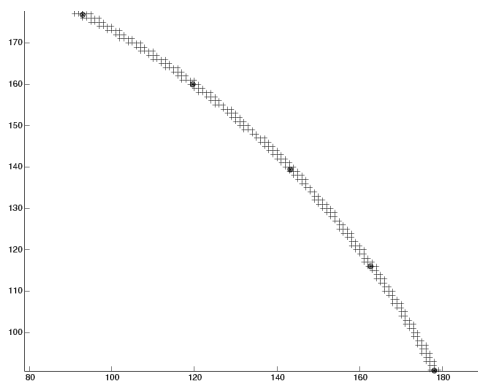}
\vspace{-0,5cm}
\caption{Valid representatives of $\mathbb X_1$}\label{ring}
\end{center}
\end{figure}

The computational timings can be drastically reduced if we perform
a grid procedure before applying AA or DA (see Section~\ref{grid_sect}).
Let us consider two cases where computation time was high: AA with
$\varepsilon=64$, and DA with $\varepsilon=2$.  In the case AA with
$\varepsilon=64$, we make a first reduction of the data using a grid whose
balls have a weighted radius of $1/4$; the computation takes $0.14$ seconds and
produces~$48$~points.  Now AA is applied to this result, and produces an
output of~$13$~points in $0.01$ seconds~---~overall far faster than applying AA
directly.  However, the final result is less accurate than that obtained by
applying AA directly. \\
The same remarks hold for the test with DA and $\varepsilon=2$: using a grid
whose balls have a weighted radius of $1/2$ we obtain $1657$ points in $0.2$
seconds; then the execution of DA on this output takes $83$ seconds to
return $466$ points.  Once again, a drastic reduction in time at the cost
of a lower quality result.
\end{ex}

\goodbreak
\begin{ex}\label{ex_zip}
{\bf Example of the  ``zip"}\\
This first example illustrates the necessity of the test at Step~AA5 of
AA. Indeed, if the condition is not checked the algorithm builds
a partition consisting of not collapsable sets.\\
Let $\mathbb X^\varepsilon$ be a set of empirical points 
whose tolerance is $\varepsilon = (2.199,2.199)$ and whose set of specified 
values $\mathbb X \subseteq \R^2$ is given by:
$$
\mathbb X = \{(0.1,2),\;(2,0),\;(4.2,0),\;(6.4,0),\;(8.6,0),\;(3.1,3)
\;(5.3,3),\;(7.5,3)\}
$$
Applying AA to the set $\mathbb X^\varepsilon$ we obtain the following partition of $\mathbb X$
\begin{eqnarray*}
 \bigl\{\{(0.1,2), (3.1,3)\},\{(2,0),(4.2,0)\},\{(6.4,0),(8.6,0)\},\{(5.3,3),(7.5,3)\}\bigr\}
\end{eqnarray*}
for which the set of specified values of the valid representatives is
$$ 
\mathbb Y=\{(1.6,2.5),\;(3.1,0),\; (7.5,0),\;(6.4,3)\}
$$
However, if we check only the distance between the centroids in step~AA5, all the elements of 
$\mathbb X^\varepsilon$ are placed in a single set which is obviously not collapsable.
\end{ex}

\bigskip
\begin{ex}\label{three}
{\bf Example of the ``three-pointed star"}\\
We have seen that AA always produces a partition into collapsable sets
such that no pair can be unified into a collapsable set.  In most
cases the partition produced by DA also enjoys this property; however,
this is not true in general.  Such a situation is shown in this example.\\
Let $\mathbb X^\varepsilon$ be a set of $6$ empirical points 
whose tolerance is $\varepsilon = (1,1)$ and whose set of specified values $\mathbb X \subseteq \R^2$
is given by:
$$
\mathbb X=\{(0.577,0.99),\;(0.577,-0.99),\;(0,0.0001),\;(0,0),\;(-1.1551,0),
\;(-1.155,0)\}
$$
Applying both AA and DA we obtain the two different partitions 
$\mathcal L_A$ and $\mathcal L_D$:
\begin{eqnarray*}
\mathcal L_A &=&\bigl\{\{(0.577,-0.99)\},\\
& &\phantom{\bigl\{}\{(0.577,0.99),(0,0.0001),(0,0)\},\\
& &\phantom{\bigl\{} \{(-1.1551,0),(-1.155,0)\}\bigr\} \\
\mathcal L_D &=&\bigl\{  \{(0.577,-0.99)\},\\
& &\phantom{\bigl\{}\{(0.577,0.99)\},\\
& &\phantom{\bigl\{} \{(0,0.0001),(0,0),(-1.1551,0),(-1.155,0)\}\bigr\}
\end{eqnarray*}
associated to the valid representatives whose specified values are
\begin{eqnarray*}
\mathbb Y_A&=&\{(0.577,-0.99), (0.192333,0.330033),(-1.15505,0)\}\\
\mathbb Y_D&=&\{(0.577,0.99),(0.577,-0.99), (-0.577525,0.000025)\}
\end{eqnarray*}
respectively.  It is trivial to verify that the elements of $\mathcal L_A^\varepsilon$ 
are pairwise not unifiable into a collapsable set, while 
the same property does not hold for the partition $\mathcal L_D^\varepsilon$ since 
$\{ (0.577,-0.99)^\varepsilon\}\cup\{(0.577,0.99)^\varepsilon \}$ is a collapsable set. 
\end{ex}

%% file: conclusions.tex
\section{Conclusions}

In this paper a new approach to reducing redundancy in sets of noisy
data is described.  The key idea is to work with empirical points,
\ie~taking into consideration the componentwise tolerances on the
input data.  The two algorithms presented are included in CoCoALib
which is available from the web site~\cite{Co}.

The experimental results points out that it is faster to use DA when
the set of empirical data is dense enough, since only a few splittings
of the original set are needed.  Conversely, when the points are
well-separated, AA is preferable, as the final partition consists of a
large number of sets and the algorithm will perform few iterations.
The very quick grid method can be used to estimate the number of final
partitions, and thus guide the choice between AA and DA.